\documentclass[10pt,reqno]{amsart}

\usepackage[utf8]{inputenc}
\usepackage{comment}
\usepackage{amsmath}
\usepackage{mathtools}
\usepackage{graphicx}
\usepackage{amsfonts}
\usepackage{amsthm}
\usepackage{amssymb}
\usepackage{mathrsfs}	
\numberwithin{equation}{section}
\usepackage[all]{xy}
\usepackage{color}
\usepackage[left=3cm,right=3cm]{geometry}
\usepackage{bm} 
\usepackage{wasysym}
\usepackage{ytableau}
\usepackage{cleveref}

\makeatletter
\newcommand{\extp}{\@ifnextchar^\@extp{\@extp^{\,}}}
\def\@extp^#1{\mathop{\bigwedge\nolimits^{\!#1}}}
\makeatother

\DeclareMathOperator{\ch}{ch}

\DeclareMathOperator{\GL}{GL}
\DeclareMathOperator{\Der}{Der}
\DeclareMathOperator{\oH}{H}
\DeclareMathOperator{\Hook}{Hook}
\DeclareMathOperator{\hook}{h}
\DeclareMathOperator{\rk}{rk}

\newcommand{\aaa}{\mathbf{a}}
\newcommand{\bbb}{\mathbf{b}}
\newcommand{\ccc}{\mathbf{c}}
\newcommand{\C}{\mathsf{C}}

\newcommand{\D}{\mathsf{D}}
\newcommand{\E}{\mathsf{\mathbf{E}}}

\newcommand{\power}{\mathsf{\mathbf{p}}}
\newcommand{\q}{\mathsf{q}}

\newcommand{\Schur}{\mathsf{\mathbf{S}}}
\newcommand{\schur}{\mathsf{\mathbf{s}}}

\newcommand{\Q}{\mathbb{Q}}

\newcommand{\N}{\mathbb{N}}

\newtheorem{lemma}{Lemma}[section]
\newtheorem{proposition}[lemma]{Proposition}
\newtheorem{thm}[lemma]{Theorem}
\newtheorem*{thmA}{Theorem A}
\newtheorem*{thmB}{Theorem B}
\newtheorem{corollary}[lemma]{Corollary}

\theoremstyle{definition}
\newtheorem{definition}[lemma]{Definition}
\newtheorem{example}[lemma]{Example}
\newtheorem{remark}[lemma]{Remark}

\title{Chern character of Schur bundles}

\author{Alessandro D'Andrea}

\author{Enrico Fatighenti}

\author{Claudio Onorati}

\address{\newline Alma Mater studiorum Università di Bologna, Dipartimento di Matematica, Piazza di Porta San Donato 5, 40126 Bologna, Italy} \email[A.~D'Andrea]{a.dandrea@unibo.it} \email[E.~Fatighenti]{enrico.fatighenti@unibo.it} \email[C.~Onorati]{claudio.onorati@unibo.it}
\date{}

\begin{document}
	
	\begin{abstract}
		We provide explicit formulas for computing the Chern character of Schur bundles $\Schur^\alpha E$ in terms of that of $E$. 
	\end{abstract}
	
	\maketitle
	
	\tableofcontents

	\section*{Introduction}
	
	Given a complex holomorphic vector bundle $E$ on a compact complex manifold $X$, one can associate to it a \emph{Schur bundle}
	\[ \Schur^\alpha E, \]
	where $\alpha$ is a Young diagram of length at most $\rk(E)$.
	If we denote by $V$ any fibre of $E$, then the corresponding fibre of $\Schur^\alpha E$ is the irreducible Schur $\GL(V)$-module $\Schur^\alpha V$.
	
	An important problem in computational algebraic geometry is to compute Chern classes of $\Schur^\alpha E$ in terms of those of $E$. It is in principle possible to perform such computations using the splitting principle, but the combinatorics behind this approach is very cumbersome, see for example \cite{Dragutin, Manivel}, where the cases of symmetric and alternating powers were treated.
	In this short paper we provide a compact and effective solution to this problem.
	
	\begin{thmA}[Theorem~\ref{thm:Chern}]
		For any $n\geq0$, the $n^{\operatorname{th}}$ Chern character coefficient of $\Schur^\alpha E$ is
		\[ \ch_n(\Schur^\alpha E)=\sum_{\nu\vdash n} \rk(\D^\nu(\schur_\alpha))\ch^{(\nu)}(E). \]
	\end{thmA}
	Let us explain the formula above:
	\begin{itemize}
		\item the sum runs over all partitions $\nu=(1^{\nu_1},2^{\nu_2},\dots)$ of the number $n$;
		\item $\schur_\alpha$ is the Schur polynomial in $\rk(E)$ variables associated to the Young diagram $\alpha$;
		\item $\D^\nu(\schur_\alpha)$ is the $\nu$-\emph{derived symmetric polynomial} which we define in Section~\ref{section:D};
		\item $\rk\colon\Q[x_1,\dots,x_r]^{\mathfrak{S}_r}\to\Q$ is the function associating to each Schur polynomial the dimension of the corresponding irreducible polynomial $\GL(r)$-representation;
		\item $\ch^{(\nu)}(E)=\frac{\ch_1(E)^{\nu_1}}{\nu_1!}\frac{\ch_2(E)^{\nu_2}}{\nu_2!}\cdots$ are the divided powers of the Chern character coefficients.
	\end{itemize}
	
	The quantity $\D^{\nu}(\schur_\alpha)$, which is the new ingredient in the formula, is a purely combinatorial object which is easily computed in explicit examples, see Examples~\ref{example:2} and \ref{example:3} and Proposition~\ref{prop:ch of sym}.
	\medskip
	
	The proof of Theorem~\ref{thm:Chern} is simple and only uses standard properties of the Chern character, such as multiplicativity and additivity. To emphasise this aspect, thus making the proof more transparent, we undress it of all the superficial structures and reduce it to a statement about polynomials in countably many variables. This is the content of our key Lemma~\ref{key lemma}, from which all the results in the paper are deduced.
	\medskip
	
	Notice that this abstract point of view is not only an {\em exercice de style}. The same proof in fact ensures that the result holds in any category where suitable characteristic classes can be defined. For example, it holds for proper algebraic varieties and locally free sheaves over them. Similarly, it also holds for characteristic classes valued in a Chow ring, instead of a cohomology ring. 
	As an example of characteristic classes that differs from the Chern character, we also prove the analogous statement for the \emph{twisted Chern character} 
	\[ k(F)=\exp\left(-\frac{c_1(F)}{\rk(F)}\right)\ch(F) \]
	considered in \cite{Markman} as characteristic classes of projective bundles.
	
	\begin{thmB}[Theorem~\ref{thm:k}]
		For any $n\geq0$, the $n^{\operatorname{th}}$ twisted Chern character coefficient of $\Schur^\alpha E$ is
		\[ k_n(\Schur^\alpha E)=\sum_{\nu\vdash n}\rk(\D^\nu(\schur_\alpha))\,k^{(\nu)}(E), \]
		where the sum may be taken over the partitions of $n$ into parts that are all strictly larger than $1$.
	\end{thmB}
	
	In our previous work \cite{DAFOa} we proved formulas expressing $k_2(\Schur^\alpha E)$ and $k_3(\Schur^\alpha E)$ as certain multiples of $k_2(E), k_3(E)$. Theorem B reinterprets those results as 
	\begin{align*}
		k_2(\Schur^\alpha E) = & \, \rk(\D_2(\schur_\alpha)) \, k_2(E)\, , \\
		k_3(\Schur^\alpha E) = & \, \rk(\D_3(\schur_\alpha)) \, k_3(E)\, ,
	\end{align*}
	and provides generalizations to all twisted Chern character coefficients $k_n$, see \S~\ref{section:k}.
	
	As already mentioned, concrete formulas for the Chern character of symmetric and exterior powers were found by Svrtan in \cite{Dragutin} and Manivel in \cite[Section~2]{Manivel}. We recover those formulas in Proposition~\ref{prop:ch of sym} and Proposition~\ref{prop:ch of wedge}, and extend them to twisted Chern characters.
	
	The formula for the first Chern class of a Schur bundle has long been folklore knowledge, and its first rigorous proof can be found in \cite{Rubei}. We present in Corollary~\ref{corollary:Rubei} an alternative short proof. 
	
	\subsection*{Plan of the paper}
	The technical core of the paper is contained in \S~\ref{section:key lemma}, where the Key Lemma is proved. In \S~\ref{section:symmetric functions} we recover some important facts about the ring of symmetric functions, that will be used in \S~\ref{section:applications}, where Theorems A and B are proved. Finally, in \S~\ref{section:Sym and Wedge} we prove formulas for Chern character and twisted Chern character of symmetric and exterior powers.
	
	\subsection*{Acknowledgements}
	We are glad to thank Laurent Manivel for comments on a first draft of the paper.
	This research has been partially funded by the European Union - NextGenerationEU under the
	National Recovery and Resilience Plan (PNRR) - Mission 4 Education and research - Component 2
	From research to business - Investment 1.1 Notice Prin 2022 - DD N. 104 del 2/2/2022, from title
	“Symplectic varieties: their interplay with Fano manifolds and derived categories”, proposal code
	2022PEKYBJ – CUP J53D23003840006.
	The second and third authors are members of the INDAM-GNSAGA group.
	
	\section{The Key Lemma}\label{section:key lemma}
	Let $R = \Q[x_1,x_2,\dots]$ be a polynomial ring {\em with an explicit choice} of countably many free commutative $\Q$-algebra generators $x_k$, $k\geq1$; we think of $R = \bigoplus R_d$ as a graded ring in which each indeterminate $x_i$ is homogeneous of any arbitrarily chosen degree.\footnote{The only role played by the grading of $R$ is in extending our statements below from monomials to all homogeneous polynomials, so the nature of the grading or the monoid/group used to grade $R$ are not really relevant.}
	We are also given a commutative unital $\N$-graded $\Q$-algebra $A=\bigoplus_{n\geq0}A_n$, where $\Q \subset A_0$, whose completion we denote by $\widehat{A}=\prod_{n\geq0}A_n$. 
	Finally, we fix a collection of positive rational numbers $\q=\{q_{n,k}\}_{n\geq1,\, k\geq1}$ with $q_{n,1}=1$.
	
	\begin{definition}\label{defn:char class}
		Consider a map $\C: R = \Q[x_1, x_2, \dots] \to \widehat{A}$ whose projection to each graded component we denote by $\C_n\colon R \to A_n$, so that
		$$\C = \sum_{n \geq 0} \C_n.$$
		We say that $\C$ is an $\widehat{A}$-\emph{valued} $\q$-\emph{characteristic class} if
		\begin{enumerate}
			\item $\C_0$ is a ring homomorphism;
			\item $\C$ is multiplicative, i.e. $\C(fg) = \C(f)\C(g)$ for all $f, g \in R$;
			\item the restriction of $\C$ to each graded component $R_d$ is $\Q$-linear;
			\item for every $n\geq1$ and every $k\geq1$, we have
			\[ \C_n(x_k)=q_{n,k}\C_n(x_1). \]
		\end{enumerate}
	\end{definition}
	
	\begin{remark}\label{rmk:same-degree additivity}
		We refer to property (3) above as \emph{same-degree additivity}. In our main application $\C$ will be the Chern character, which is additive; nevertheless, allowing this more general setting gives us the possibility to extend the results to other classes, see for example \S~\ref{section:k}.
	\end{remark}
	
	For each $n\geq1$, let us now consider the derivations
	\begin{equation}\label{eqn:der with q}
		\delta_{\mathsf{q},n}=\sum_{k\geq1}q_{n,k}\frac{\partial}{\partial x_k}\in \Der R.
	\end{equation}
	Notice that each $\delta_{\mathsf{q},n}$ is locally nilpotent and well defined, since for each polynomial $f\in\Q[x_1,x_2,\dots]$ only finitely many of the $\frac{\partial f}{\partial x_k}$ are non-zero.
	
	The formal sum $\sum_{n \geq 1} \C_n(x_1) \cdot \delta_{\mathsf{q},n}$ is then a (locally nilpotent, converging) derivation of $R_{\widehat A} = \widehat A \otimes_\Q R = \widehat A[x_1, x_2, \dots]$ and its exponential is a well defined automorphism of $R_{\widehat A}$, whose restriction to $R \subset R_{\widehat A}$ yields a $\Q$-algebra homomorphism
	\[ \exp\left(\sum_{n\geq1}\C_n(x_1)\cdot\delta_{\mathsf{q},n}\right)\colon\Q[x_1,x_2,\dots]\longrightarrow R_{\widehat A} = \widehat{A}\otimes_\Q \Q[x_1,x_2,\dots]. \]
	Finally, we abuse notation and also denote by $\C_0$ the composition
	$$\widehat A \otimes_\Q \Q[x_1,x_2,\dots] \stackrel{1 \otimes \C_0}{\longrightarrow} \widehat A \otimes_\Q A_0 \stackrel{\mu}{\longrightarrow} \widehat A,$$
	where $\mu$ is the multiplication map of $\widehat A$. Notice that since $\widehat A$ is commutative, the multiplication map $\mu$ is a ring homomorphism; consequently, the composition
	\[ \C_0\circ\exp\left(\sum_{n\geq1}\C_n(x_1)\cdot\delta_{\mathsf{q},n}\right)\colon R\longrightarrow \widehat{A}  \]
	is a well-defined $\Q$-algebra homomorphism by construction.
	
	\begin{lemma}[Key Lemma]\label{key lemma}
		Let $\C$ be an $\widehat{A}$-valued $\q$-characteristic class. Then equality
		\begin{equation}\label{key} \C = \C_0\circ\exp\left(\sum_{n\geq1}\C_n(x_1)\cdot\delta_{\mathsf{q},n}\right)
		\end{equation}
		holds whenever both sides are applied to any given homogeneous polynomial in $R$. If furthermore $\C$ is a $\Q$-algebra homomorphism, then equality \eqref{key} holds in general.
		\begin{proof}
			Let us prove that \eqref{key} holds when applied to the element $x_k$. The left-hand side is
			$$\C(x_k) = \sum_{n \geq 0} \C_n(x_k) = \C_0(x_k) + \sum_{n \geq 1} q_{n,k} \C_n(x_1).$$
			
			As for the right-hand side, we notice that $\delta_{\mathsf{q},n}(x_k) = q_{n,k} \in \Q$, so that $\delta_{\mathsf{q},m}(\delta_{\mathsf{q},n}(x_k)) = 0$ for all $m$ and $n$. When we apply the automorphism to $x_k$, only the constant and the linear term survive and we obtain
			$$\exp\left(\sum_{n\geq 1} \C_n(x_1) \cdot \delta_{\mathsf{q},n} \right)(x_k) = x_k + \sum_{n \geq 1} q_{n,k} \C_n(x_1).$$
			After applying $\C_0$, this coincides with the left-hand side.
			
			By multiplicativity of $\C$, Equation \eqref{key} holds for each monomial, and since the restriction of $\C$ on each graded component of $R$ is $\Q$-linear, it also holds for every homogeneous polynomial. The final claim is clear.
		\end{proof}
	\end{lemma}

	For the following corollary we use the notation $\nu=(1^{\nu_1},2^{\nu_2},\dots,h^{\nu_h})$ for partitions. In particular, $\nu\vdash n$ means that $\sum_i i\nu_i=n$.
	Then we set
	\begin{equation}\label{eqn: D mu} \delta_{\mathsf{q}}^\nu=\delta_{\mathsf{q},1}^{\nu_1}\circ\cdots\circ\delta_{\mathsf{q},h}^{\nu_h}.  \end{equation}
	Notice that, for any $n,m\in\mathbb{Z}$, the operators $\delta_{\mathsf{q},n}$ and $\delta_{\mathsf{q},m}$ commute, i.e.\ $\delta_{\mathsf{q},n}\delta_{\mathsf{q},m}=\delta_{\mathsf{q},m}\delta_{\mathsf{q},n}$.
	
	\begin{corollary}\label{cor:key lemma}
		Let $\C$ be a $\hat{A}$-valued $\q$-characteristic class. Then for every $n\geq0$ and every homogeneous $f\in\Q[x_1,x_2,\dots]$
		\[ \C_n(f)=\sum_{\nu\vdash n}\C_0(\delta_{\mathsf{q}}^\nu(f))\C^{(\nu)}(x_1), \]
		where $\C^{(\nu)}=\prod_{i}\frac{\C_i^{\nu_i}}{\nu_i!}$ denotes divided powers of coefficients of $\C$. \qed
	\end{corollary}

	\section{Symmetric functions}\label{section:symmetric functions}
	In this short section we recall standard facts about symmetric functions, see \cite{Macdonald} for a reference.
	
	\subsection{}
	The following are the notations we use in this and the next section.
	\begin{itemize}
		\item Young diagrams will be denoted by the letter $\alpha$. We also represent Young diagrams via vectors whose entries denote the length of rows of $\alpha$. The \emph{size} $|\alpha|$ of a Young diagram is the number of boxes or, equivalently, the sum of the entries in the vector representation. The \emph{length} of $\alpha$ is the number of rows or, equivalently, the number of non-zero entries in its vector representation.
		\item The notation $(i,j)\in\alpha$ means that we are considering the box in position $(i,j)$ in $\alpha$. The upper-left box has coordinates $(1,1)$ and we increase each coordinate by one at every step downwards or rightwards, accordingly.
		\item Given $(i,j)\in\alpha$, we denote by $\Hook_{(i,j)}(\alpha)$ the \emph{hook} centred in position $(i,j)$ and by $\hook_{(i,j)}(\alpha)$ the corresponding \emph{hook number}, i.e.\ its size.
		\item Partitions will be denoted by $\nu=(1^{\nu_1},2^{\nu_2},\dots)$ and their size is $|\nu|=\sum_i i\nu_i$.
		\item $\Lambda=\bigoplus\Lambda^d$ is the graded ring of symmetric functions with coefficients in $\Q$.
		\item Schur functions are denoted by $\schur_\alpha$ and are indexed by Young diagrams. For $|\alpha|=d$ we have $\schur_\alpha\in\Lambda^d$.
		\item Power-sum functions are $\power_k(x_1,x_2,\dots)=\sum_{i\geq1}x_i^k\in\Lambda^k$.
	\end{itemize}
	
	\subsection{}
	
	Let $\Lambda=\bigoplus\Lambda^d$ be the $\N$-graded ring of symmetric functions with coefficients in $\Q$. Two sets of generators of $\Lambda$ are given by \emph{Schur functions} $\schur_\alpha$ and \emph{power-sum functions} $\power_k$. A Schur function $\schur_\alpha$ is homogeneous of degree $|\alpha|$, while a power-sum function $\power_k$ is homogeneous of degree $k$.
	
	Both Schur functions $\schur_\alpha$ with $|\alpha|=d$ and monomials in power-sums $\power^\mu$ with $|\mu|=d$ constitute $\Q$-linear bases of each homogeneous component $\Lambda^d$. In particular power-sums are algebraically independent generators of 
	\[ \Lambda=\Q[\power_1,\power_2,\dots], \]
	so that $\Lambda$ may be thought of as a ring of polynomials with rational coefficients in countably many variables $\power_1, \power_2, \dots$
	
	\subsection{}\label{section:D}
	
	For any integer $n$, define the derivations
	\begin{equation}\label{eqn:D}
		\D_n\colon \Lambda\to\Lambda\,, \qquad \D_n=\sum_{k\geq0}k^n\frac{\partial}{\partial\power_k}.
	\end{equation}
	In \cite{DAFOb} we introduced an analogous derivation $\D_s$, where $s$ is a formal variable, and used it 
	to prove recursive formulas that are useful in the computation of tensor products and plethysms of Schur functions. Derivations in \eqref{eqn:D} are obtained by evaluating $\D_s$ at $s=n \in \N$.
	
	\begin{remark}\label{rmk:D original}
		Notice that the $\D_n$ are the derivations $\delta_{\mathsf{q},n}$ we obtain from (\ref{eqn:der with q}), when choosing $q_{n,k} = k^n$.
	\end{remark}
	
	Next, we want to explain how to effectively compute $\D_n(\schur_\alpha)$ using the following well-known combinatorial interpretation, which is explained for instance in \cite[Section~3]{DAFOb} (see also \cite[Proposition~2.10]{DAFOb}).
	For each $(i,j)\in\alpha$ we denote by $\alpha_{(i,j)}$ the Young diagram obtained by $\alpha$ by removing the border-strip having the same ends as $\Hook_{(i,j)}(\alpha)$.\footnote{Equivalently, we may delete $\Hook_{(i,j)}(\alpha)$ and move the lower component up-left if removing the hook disconnects the diagram.}
	
	\begin{example}\label{example:(3,2,2)}
		If $\alpha=(3,2,2)$, then $\alpha_{(1,2)}=(3,1)$.
	\end{example}
	
	\begin{proposition}[\protect{\cite[(3.1)]{DAFOb}}]\label{prop:D geometrica}
		The following equality holds,
		\[ \D_n(\schur_\alpha)=\sum_{(i,j)\in\alpha} (-1)^{\varepsilon(i,j)}(\hook_{(i,j)})^{n-1}\cdot \schur_{\alpha_{(i,j)}}, \]
		where $\varepsilon(i,j)$ is the number of boxes strictly below $(i,j)$.
	\end{proposition}
	
	\begin{example}
		Continuing Example~\ref{example:(3,2,2)}, we have
		\[ \D_n(\schur_{(3,2,2)})=\schur_{(3,2,1)}+\schur_{2,2,2}+(\schur_{(3,2)}-\schur_{(3,1,1)})2^{n-1}-\schur_{3,1}3^{n-1}+\schur_{1,1,1}4^{n-1}+\schur_{(1,1)}5^{n-1}\,. \]
	\end{example}
	
	\begin{example}\label{example:D of Sym and Wedge}
		For symmetric and exterior powers, we have
		\[ \D_n(\schur_m)=\sum_{k=1}^mk^{n-1}\schur_{m-k}\qquad\mbox{and}\qquad \D_n(\schur_{1^m})=\sum_{k=1}^m(-1)^{k+1}k^{n-1}\schur_{1^{m-k}}. \]
	\end{example}
	
	As in (\ref{eqn: D mu}), we put $\D^\mu=\D_\ell^{\mu_{\ell}}\cdots\D_1^{\mu_1}$, where $\mu=(1^{\mu_1},\dots,\ell^{\mu_\ell})$ is a partition, and notice that $\D_n\cdot\D_m=\D_m\cdot\D_n$ for every $m$ and $n$.
	
	\begin{example}\label{example:D 1 k}
		Since $\D_1(\schur_m)=\schur_{m-1}+\schur_{m-2}+\cdots+1$, it is easy to see that, for every $k>0$,
		\[ \D_1^k(\schur_m)=\sum_{j=0}^{m-k}\binom{m-1-j}{k-1}\schur_j. \]

		Similarly, it is also easy to see that $\D_2^k(\schur_m)=\D_1^{2k}(\schur_{m+k})$.
	\end{example}
	
	\begin{remark}
		If $|\alpha|=n$, then $\D^\mu(\schur_\alpha)=0$ for every $\mu=(1^{\mu_1},\dots,\ell^{\mu_\ell})$ such that $\mu_1+\cdots+\mu_\ell>n$. (If we express $\mu$ as a Young diagram, then $\mu_1+\cdots+\mu_\ell$ is the number of rows of the diagram.)
	\end{remark}
	
	\section{Chern characters of Schur bundles}\label{section:applications}
	
	Throughout this section we assume that $X$ is a connected compact complex manifold and choose the even cohomology ring $\oH^{2\ast}(X,\Q)$ as the graded commutative algebra $A$ to use in Lemma \ref{key lemma}.
	Notice that since $\oH^{2n}(X,\Q)$ is nonzero for finitely many $n$, then $A$ equals its completion $\widehat A$. 
	
	\subsection{Chern character}\label{section:Chern}
	Let $E$ be a rank $r$ holomorphic vector bundle on $X$.
	Consider the following homomorphism
	\begin{equation}\label{eqn:Chern} 
		\C^E\colon \Lambda = \Q[\power_1, \power_2, \dots] \to \oH^{2\ast}(X,\Q),\qquad \schur_\alpha\mapsto\ch(\Schur^\alpha E). 
	\end{equation}
	
	Let us check that $\C^E$ is a {\em characteristic class} as from Definition~\ref{defn:char class}. 
	
	\begin{lemma}
		Set $\q=\{k^n\}_{n\geq1,k\geq1}$. Then $\C^E$ is a $\oH^{2\ast}(X,\Q)$-valued $\q$-characteristic class.
		\begin{proof}
			The first three points of Definition~\ref{defn:char class} are straightforward; notice that in this case $\C^E$ is fully additive, not just same-degree additive. Finally, the equality $\C_n^E(\power_k)=k^n\C_n^E(\power_1)$ can be checked by hand, or one can look at \cite[Example~1.5]{DAFOa}.
		\end{proof}
	\end{lemma}
	
	\begin{remark}
		The set $\q=\{k^n\}_{n\geq1,k\geq1}$ is the same used to define the derivation $\D_n$ in \S~\ref{section:D}, see Remark~\ref{rmk:D original}.
	\end{remark}
	
	\begin{remark}
		The proof of (4) in Definition~\ref{defn:char class} makes implicit use of the splitting principle: this is what we need in order to apply \cite[Example~1.5]{DAFOa}. 
		We stress the fact that this is the only point where the splitting principle is used.
	\end{remark}
	
	Recall that, by the hook length dimension formula, one has
	\[ \C^E_0(\schur_\alpha)=\rk(\Schur^\alpha E)=\prod_{(i,j)\in\alpha}\frac{r+i-j}{\hook_{(i,j)}(\alpha)}, \]
	so that, for instance,
	\[ 
	\rk(\Schur^{(2,1)} E)= \frac{r}{3}\cdot\frac{r+1}{1}\cdot\frac{r-1}{1}=2\binom{r+1}{3}\quad\mbox{ and }\quad \rk(\Schur^{(3,2)}E)=\frac{(r+2)(r+1)r^2(r-1)}{24}.
	\]
	
	As a direct application of Corollary~\ref{cor:key lemma}, we get the following formula for the Chern character of Schur bundles.
	
	\begin{thm}\label{thm:Chern}
		With notations as above, we have
		\[ \ch_n(\Schur^\alpha E)=\sum_{\nu\vdash n}\rk(\D^\nu(\schur_\alpha))\,\ch^{(\nu)}(E). \]
	\end{thm}
	
	Let us see some consequences and explicit examples.
	
	\begin{corollary}[\cite{Rubei}]\label{corollary:Rubei}
		With notations as above,
		\[ \ch_1(\Schur^\alpha E)=|\alpha|\frac{r_\alpha}{r}\ch_1(E), \]
		where $r_\alpha=\rk(\Schur^\alpha E)$.
		\begin{proof}
			From Theorem~\ref{thm:Chern} we obtain that
			$$\ch_1(\Schur^\alpha E) = \rk(\D_1(\schur_\alpha)) \ch_1(E).$$
			As $$\D_1 = \sum_{k \geq 0} \frac{\partial}{\partial \power_k}$$ and $\rk(\power_k) = r$ for every $k$, 
			then
			$$\rk(\D_1(\schur_\alpha)) = \frac{1}{r} \rk\left(\sum_k \power_k \frac{\partial}{\partial \power_k} (\schur_\alpha)\right).$$ 
			The claim follows now by noticing that $\sum_k \power_k {\partial}/{\partial \power_k}$  
			acts on homogeneous elements via multiplication by their degree.
		\end{proof}
	\end{corollary}
		
	We wish to point out the effectiveness of Theorem~\ref{thm:Chern} by explicitly computing some examples.
	
	\begin{example}[Schur bundles of size $2$]\label{example:2}
		Let us compute $\ch(\Schur^{\alpha}E)$ when $|\alpha|=2$. We have only two cases, namely the symmetric power $\Schur^2E$ and the exterior power $\extp^2 E$, corresponding to $\alpha=(1,1)$. 
		\begin{description}
			\item[$\alpha=(2)$] it is easy to see that the only non-zero derivatives of $\schur_2$ are
			\[ \D_a(\schur_2)=\schur_1+2^{a-1}\implies \rk(\D_a(\schur_2))=r+2^{a-1} \qquad\mbox{ and }\qquad\D_a\D_b(\schur_2)=1. \]
			Then the first values of $\ch_n(\Schur^2 E)$ are:
			\begin{align*}
				\ch_1(\Schur^{2}E) = & \, (r+1)\ch_1(E), \\
				\ch_2(\Schur^{2}E) = & \, (r+2)\ch_2(E)+\frac{\ch_1(E)^2}{2}, \\
				\ch_3(\Schur^{2}E) = & \, (r+4)\ch_3(E)+\ch_2(E)\ch_1(E), \\
				\ch_4(\Schur^{2}E) = & \, (r+8)\ch_4(E)+\ch_3(E)\ch_1(E)+\frac{\ch_2(E)^2}{2}, \dots
			\end{align*}
			
			\item[$\alpha=(1,1)$] it is easy to see that the only non-zero derivatives of $\schur_{(1,1)}$ are
			\[ \D_a(\schur_{(1,1)})=\schur_1-2^{a-1}\implies \rk(\D_a(\schur_2))=r-2^{a-1}\qquad\mbox{ and }\qquad\D_a\D_b(\schur_{(1,1)})=1. \]
			Then the first values of $\ch_n(\extp^2 E)$ are:
			\begin{align*}
				\ch_1(\extp^2E) = & \, (r-1)\ch_1(E), \\
				\ch_2(\extp^2E) = & \, (r-2)\ch_2(E)+\frac{\ch_1(E)^2}{2}, \\
				\ch_3(\extp^2E) = & \, (r-4)\ch_3(E)+\ch_2(E)\ch_1(E), \\
				\ch_4(\extp^2E) = & \, (r-8)\ch_4(E)+\ch_3(E)\ch_1(E)+\frac{\ch_2(E)^2}{2}, \dots
			\end{align*}
		\end{description}
	\end{example}
	
	\begin{example}[Schur bundles of size $3$]\label{example:3}
		In this case we have three cases, namely $\alpha=(3),(2,1)$ and $(1,1,1)$. 
		\begin{description}
			\item[$\bullet\;\alpha=(3)$] It is easy to compute
			\[ \D_a(\schur_3)=\schur_2+\schur_1 2^{a-1}+3^{a-1},\quad \D_a\D_b(\schur_3)=\schur_1+2^{a-1}+2^{b-1}\quad\mbox{and}\quad \D_a\D_b\D_c(\schur_3)=1, \]
			from which we deduce the first values of $\ch_n(\Schur^3E)$:
			\begin{align*}
				\ch_1(\Schur^{3}E) = & \, \binom{r+2}{2}\ch_1(E), \\
				\ch_2(\Schur^{3}E) = & \, \binom{r+3}{2}\ch_2(E)+(r+2)\frac{\ch_1(E)^2}{2}, \\
				\ch_3(\Schur^{3}E) = & \, \frac{(r+3)(r+6)}{2}\ch_3(E)+(r+3)\ch_2(E)\ch_1(E) + \frac{\ch_1(E)^3}{6}, \\
				\ch_4(\Schur^{3}E) = & \, \frac{r^2+17r+54}{2}\ch_4(E)+(r+5)\ch_3(E)\ch_1(E)+(r+4)\frac{\ch_2(E)^2}{2}+\frac{\ch_2(E)\ch_1(E)^2}{2}.
			\end{align*}
			
			\item[$\bullet\;\alpha=(2,1)$] It is easy to see that
			\[ \D_{a}(\schur_{(2,1)})=\schur_{2}+\schur_{(1,1)}-3^{a-1},\quad \D_a\D_b(\schur_{(2,1)})=2\schur_1\qquad \mbox{and}\qquad \D_a\D_b\D_c(\schur_{(2,1)})=2, \]
			from which we deduce the first values of $\ch_n(\Schur^{(2,1)}E)$:
			\begin{align*}
				\ch_1(\Schur^{(2,1)}E) = & \, (r^2-1)\ch_1(E), \\
				\ch_2(\Schur^{(2,1)}E) = & \, (r^2-3)\ch_2(E)+2r\frac{\ch_1(E)^2}{2}, \\
				\ch_3(\Schur^{(2,1)}E) = & \, (r^2-9)\ch_3(E)+2r\ch_2(E)\ch_1(E) + 2 \frac{\ch_1(E)^3}{6}, \\
				\ch_4(\Schur^{(2,1)}E) = & \, (r^2-27)\ch_4(E)+2r\ch_3(E)\ch_1(E)+2r\frac{\ch_2(E)^2}{2}+2\frac{\ch_2(E)\ch_1(E)^2}{2}, \dots
			\end{align*}
			
			\item[$\bullet\;\alpha=(1,1,1)$] It is easy to compute
			\[ \D_a(\schur_{(1,1,1)})=\schur_{(1,1)}-\schur_1 2^{a-1}+3^{a-1},\quad \D_a\D_b(\schur_3)=\schur_1-2^{a-1}-2^{b-1}\quad\mbox{and}\quad \D_a\D_b\D_c(\schur_3)=1, \]
			from which we deduce the first values of $\ch_n(\extp^3E)$:
			\begin{align*}
				\ch_1(\extp^3E) = & \, \binom{r-1}{2}\ch_1(E), \\
				\ch_2(\extp^3E) = & \, \binom{r-2}{2}\ch_2(E)+(r-2)\frac{\ch_1(E)^2}{2}, \\
				\ch_3(\extp^3E) = & \, \frac{(r-3)(r-6)}{2}\ch_3(E)+(r-3)\ch_2(E)\ch_1(E) + \frac{\ch_1(E)^3}{6}, \\
				\ch_4(\extp^3E) = & \, \frac{r^2-17r+54}{2}\ch_4(E)+(r-5)\ch_3(E)\ch_1(E)+(r-4)\frac{\ch_2(E)^2}{2}+\frac{\ch_2(E)\ch_1(E)^2}{2}, \dots
			\end{align*}
		\end{description}
	\end{example}
	
	\subsection{Twisted Chern character}\label{section:k}
	
	In \cite{Markman}, Markman introduces the \emph{twisted Chern character} of a vector bundle $F$ as 
	\[ k(F)=\exp\left(-\frac{\ch_1(F)}{\rk(F)}\right)\ch(F). \]
	This class should be thought of as the Chern character of the projective bundle $\mathbb{P}E$, which motivates the terminology. For bookkeeping, let us express the first three values of $k(F)$ in terms of Chern classes:
	\[ k_0(F)=\rk(F),\quad k_1(F)=0,\quad  k_2(F)=-\frac{\ch_1(F)^2-2\rk(F)\ch_2(F)}{2\rk(F)}  \]
	and
	\[ k_3(F)=\frac{\ch_1(F)^3-3\rk(F)\ch_2(F)\ch_1(F)+3\rk(F)^2\ch_3(F)}{3\rk(F)^2}. \]
	Notice that $-2\rk(F)\, k_2(F)=\Delta(F)$ is the discriminant of $F$.
	\medskip
	
	Define
	\[ \C^{\mathbb{P}E}\colon\Lambda\to\oH^{2\ast}(X,\Q),\qquad \schur_\alpha\mapsto k(\Schur^\alpha E) = \exp\left(-\frac{|\alpha|}{r}\ch_1(E)\right)\ch(\Schur^\alpha E), \]
	where the equality follows from Corollary~\ref{corollary:Rubei}. 
	
	\begin{lemma}
		Set $\q=\{k^n\}_{n\geq1,k\geq1}$. Then $\C^{\mathbb{P}E}$ is a $\oH^{2\ast}(X,\Q)$-valued $\q$-characteristic class.
		\begin{proof}
			The first two properties of Definition~\ref{defn:char class} are straightforward. Moreover, $\C^{\mathbb{P}E}$ is same-degree additive, since the degree $d$ homogeneous component of $\Lambda$ is generated by Schur functions $\schur_{\alpha}$ with $|\alpha|=d$.
			
			Finally, let us notice that $\C^{\mathbb{P}E}_1(\power_k)=0=\C^{\mathbb{P}E}(\power_1)$, while $\C^{\mathbb{P}E}_n(\power_k)=k^n\C^{\mathbb{P}E}_n(\power_1)$, as it can be computed by hand using the same equality for $\C^E$.
		\end{proof}
	\end{lemma}
	
	As a direct application of Corollary~\ref{cor:key lemma}, we get the following formula for the class $k$ of Schur bundles. 
	
	\begin{thm}\label{thm:k}
		With notations as above,
		\[ k_n(\Schur^\alpha E)=\sum_{\nu\vdash n}\rk(\D^\nu(\schur_\alpha))\,k^{(\nu)}(E). \]
	\end{thm}
	
	\begin{example}\label{example:primi k}
		The first values of $k_n(\Schur^\alpha E)$ are easy to compute. In fact, since $k_1(E)=0$, we get
		\begin{align*}
			k_2(\Schur^\alpha E) = & \, \rk(\D_2(\schur_\alpha)) \, k_2(E)\, ; \\
			k_3(\Schur^\alpha E) = & \, \rk(\D_3(\schur_\alpha)) \, k_3(E)\, ; \\
			k_4(\Schur^\alpha E) = & \, \rk(\D_4(\schur_\alpha)) \, k_4(E) + \frac{1}{2}\rk(\D_2^2(\schur_\alpha)) k_2^2(E)\, ,
			\intertext{so that, for instance,}
			k_2(\Schur^2 E) = & \, (r+2) \, k_2(E)\, , \\
			k_3(\Schur^2 E) = & \, (r+4) \, k_3(E)\, , \\
			k_4(\Schur^2 E) = & \, (r+8) \, k_4(E) + \frac{1}{2} k_2^2(E)\, .
		\end{align*}
	\end{example}
	
	\begin{remark}
		It is worth noticing that the class $k$ is related to the logarithmic Chern character as defined in \cite{Drezet}. If we denote by $\mathsf{L}(E)$ the latter, then for any vector bundle $F$,
		\begin{equation}\label{eqn:k vs L} 
			\operatorname{Log}(k(F))=\mathsf{L}(E)-\left(\frac{\ch_1(F)}{\rk(F)}\right), 
		\end{equation}
		from which it follows, for example, that 
		\[ \mathsf{L}_2(F)=\frac{k_2(F)}{\rk(F)},\quad \mathsf{L}_3(F)=\frac{k_3(F)}{\rk(F)}\quad\mbox{and}\quad \mathsf{L}_4(F)=\frac{k_4(F)}{\rk(F)}-\frac{k_2(F)^2}{\rk(F)^2}. \]
		
		In \cite{DAFOa} we proved that the classes  $\mathsf{L}_2(F)$ and $\mathsf{L}_3(F)$ are \emph{eigen-classes} with respect to the induced action of Schur functors. Moreover, we also remarked that this claim does not extend to higher-degree classes and speculated that similar formulas should exist for the projection of the classes $\mathsf{L}_n(\Schur^\alpha E)$ onto the $\mathsf{L}_n(E)$-component. 
		
		Theorem~\ref{thm:k} solves this problem in the best possible way. In fact, as it is clear from Example~\ref{example:primi k}, both $k_2$ and $k_3$ are eigen-classes, thus recovering the result for $\mathsf{L_2}$ and $\mathsf{L_3}$ in \cite{DAFOa}. Moreover, Theorem~\ref{thm:k} provides a formula not just for the $\mathsf{L}_n(E)$-component, but for all components (up to re-ordering according to (\ref{eqn:k vs L})).
	\end{remark}
	
	\begin{remark}
		Since $\C^{\mathbb{P}E}_1(\power_k)=0=\C^{\mathbb{P}E}(\power_1)$, the numbers $q_{1,k}$ can be arbitrarily chosen. Any special choice of such numbers has the effect of modifying the single derivation $\D_1$, which is however only involved in those terms $k^{(\nu)}$ containing $k_1(E)=0$ that do, in any case, vanish altogether.     
	\end{remark}

	\section{Chern characters of symmetric and exterior powers}\label{section:Sym and Wedge}
	
	In this section we explicitly compute the (twisted) Chern character of the symmetric and exterior powers of $E$, recovering precisely the formulas obtained by Manivel in \cite{Manivel} in the symmetric case and finding expressions equivalent to those by Svrtan in \cite{Dragutin}, which use Stirling numbers of the second kind instead of Eulerian numbers. 
	
	Let us start by setting notations. If $\nu=(1^{\nu_1},\dots,\ell^{\nu_\ell})$ is a partition, then we denote by 
	\[ \aaa_{\nu}=(\overbrace{\ell,\dots,\ell}^{\nu_\ell},\dots,\overbrace{1,\dots,1}^{\nu_1}) \]
	the associated Young diagram. We denote by $a_i$ the entries of $\aaa_{\nu}$. 
	Accordingly we write
	\[ \D^\nu=\D_{\aaa_\nu}=\overbrace{\D_{\ell}\circ\cdots\circ\D_{\ell}}^{\nu_\ell}\circ\cdots\circ\overbrace{\D_{1}\circ\cdots\circ\D_{1}}^{\nu_1}. \]
	Let us set
	\[ h=h(\nu)=\nu_1+\cdots+\nu_\ell \]
	to be the height of the Young diagram $\aaa_\nu$, i.e.\ the number of rows. If $\bbb=(b_{1},\dots,b_{h})$ is an $h$-tuple of integers, we put
	\[ \overleftarrow{\bbb}=(b_{1}-1,\dots,b_{h}-1). \]
	Finally, for every $n,j\geq0$, let us denote by $\E(n,j)$ the corresponding \emph{Eulerian number}. Recall the following important property of Eulerian numbers:
	\begin{equation}\label{eqn:eulerians}
		k^a=\sum_{b=1}^{a+1}\E(a,b-1)\binom{k+b-1}{a}. 
	\end{equation}
	
	\begin{remark}\label{rmk:eulerians}
		Let us notice that $\E(n,k)=0$ for all $k\geq n$, unless $\E(0,0)$ which is set equal to $1$ by definition. Formula (\ref{eqn:eulerians}) is usually written with $b$ ranging from $1$ to $a$, but would leave the case $a=0$ out, which we need in the proof of Proposition~\ref{prop:ch of sym}.
	\end{remark}
	
	If $\bbb$ and $\ccc$ are two $h$-tuples of natural numbers, then we put
	\[ \E(\ccc,\bbb)=\prod_{i=1}^h\E(c_i,b_i). \]
	Similarly, we write 
	\[ k\leq \bbb\leq \ccc \implies k\leq b_i\leq c_i\quad \forall i=1,\dots,h\] 
	and
	\[ |\bbb|=b_1+\cdots+b_h. \]
	
	Let us state the formulas. We start with the symmetric power.
	
	\begin{proposition}\label{prop:ch of sym}
		With notations as above, we have
		\begin{equation*}
			\ch_n(\Schur^m E)=\sum_{\nu\vdash n}\sum_{1\leq\bbb\leq\aaa_\nu}\E(\overleftarrow{\aaa_\nu},\overleftarrow{\bbb})\binom{r+m+|\bbb|-\nu_1-1}{r+|\aaa_\nu|-1}\,\ch^{(\nu)}(E) 
		\end{equation*}
		and
		\begin{equation*}
			k_n(\Schur^m E)=\sum_{\substack{\nu\vdash n \\ \nu_1=0}}\sum_{1\leq\bbb\leq\overleftarrow{\aaa_\nu}}\E(\overleftarrow{\aaa_\nu},\overleftarrow{\bbb})\binom{r+m+|\bbb|-1}{r+|\aaa_\nu|-1}\,k^{(\nu)}(E). 
		\end{equation*}
	\end{proposition}
	
	For exterior products, we have instead the following statement.
	
	\begin{proposition}\label{prop:ch of wedge}
		With notations as above, we have
		\begin{equation*}
			\ch_n\left(\extp^mE\right)=\sum_{\nu\vdash n}\sum_{1\leq\bbb\leq\aaa_\nu}\E(\overleftarrow{\aaa_\nu},\overleftarrow{\bbb})\binom{r-|\aaa_\nu|}{m-|\aaa_\nu|+|\bbb|-\nu_1}\,\ch^{(\nu)}(E) 
		\end{equation*}
		and
		\begin{equation*} 
			k_n\left(\extp^mE\right)=\sum_{\substack{\nu\vdash n \\ \nu_1=0}}\sum_{1\leq\bbb\leq\overleftarrow{\aaa_\nu}}\E(\overleftarrow{\aaa_\nu},\overleftarrow{\bbb})\binom{r-|\aaa_\nu|}{m-|\aaa_\nu|+|\bbb|}\,k^{(\nu)}(E).
		\end{equation*}
	\end{proposition}
	
	\begin{remark}\label{rmk:bound}
		We can be more precise about the bound $\bbb\leq\aaa_{\nu}$. In fact, if $a_i>1$, then we have $b_i\leq a_i-1$, as it will be clear from the proof. The extremal case $b_i=a_i$ only appears when $a_i=1$, i.e.\ $a_i-1=0$, see also Remark~\ref{rmk:eulerians}. This also explains why in the formula for the twisted Chern character we can restrict to $\bbb\leq\overleftarrow{\aaa_\nu}$.
	\end{remark}
	
	\subsection{Proof of Proposition~\ref{prop:ch of sym}}\label{section:proof of sym}
	
	We will need the following well-known binomial identity.
	
	\begin{lemma}\label{lemma:binomial for sym}
		For any non-negative integers $A$, $R$ and $N$,
		\begin{equation*} 
			\sum_{i=0}^N\binom{A+i}{A}\binom{R+N-i}{R}=\binom{R+N+A+1}{R+A+1}.
		\end{equation*}
	\end{lemma}
	
	First of all, it is enough to prove that 
	\begin{equation}\label{eqn:in the proof of ch of sym} 
		\rk(\D_{\aaa_\nu}(\schur_m))=\sum_{1\leq\bbb\leq\aaa_\nu}\E(\overleftarrow{\aaa_\nu},\overleftarrow{\bbb})\binom{r+m+|\bbb|-\nu_1-1}{r+|\aaa_\nu|-1}. 
	\end{equation}
	In fact, this is exactly the coefficient in front of $\ch^{(\nu)}(E)$ in the formula for the Chern character, while for the twisted Chern character it is enough to notice that $k_1(E)=0$, so that every partition $\nu$ such that $\nu_1\geq1$ does not contribute to the formula. (The shift $\overleftarrow{\aaa_\nu}$ in the formula for the twisted Chern character is explained in Remark~\ref{rmk:bound}.)
	
	We proceed by induction on $h=h(\nu)$. As base of the induction, let us prove the case $\aaa_\nu=(a)$.
	By definition and (\ref{eqn:eulerians}),
	\[ \rk(\D_a(\schur_m))=\sum_{k=1}^m k^{a-1} \binom{r+m-k-1}{r-1}=\sum_{b=1}^{a}\E(a-1,b-1)\sum_{k=1}^m \binom{k+b-1}{a-1}\binom{r+m-k-1}{r-1}. \]
	The sum on the right can be computed using Lemma~\ref{lemma:binomial for sym}, but it is important to distinguish two cases. First, if $a=1$, then the sum reduces to the well-known equality $\sum_{k=1}^m \binom{r+m-k-1}{r-1}=\binom{r+m-1}{r}$, that agrees with (\ref{eqn:in the proof of ch of sym}), since in this case only $b=1$ is allowed and $\nu_1=1$. (Notice that in this case, contrary to the next one, the sum starts with $k=1$; this is the reason why we have an extra $-1$ in the numerator of the binomial coefficient.)
	If $a\geq 2$, then we have non-trivial summands only for $a-b\leq k\leq m$ so that we can apply Lemma~\ref{lemma:binomial for sym} with $A=a-1$, $R=r-1$ and $N=m-a+b$, yielding
	\[ \rk(\D_a(\schur_m))=\sum_{b=1}^{a-1}\E(a-1,b-1)\binom{r+m+b-1}{r+a-1}. \]
	Since $\nu_1=0$ in this case, this agrees again with (\ref{eqn:in the proof of ch of sym}). 
	
	For the inductive step, let us first assume that there exists an index $i$ such that $a_i>1$. Without loss of generality, we can assume that it is $a_1>1$. In fact, since $\D_{a_i}$ and $\D_{a_j}$ commute we can re-order the indices. Let us then write $\aaa_\nu=(a_1,\dots,a_{h})$, and let us set $\aaa'=(a_2,\dots,a_h)$. By definition we have
	\[ \rk(\D_{\aaa}(\schur_m))=\sum_{k=1}^m k^{a_{1}-1}\rk(\D_{\aaa'}(\schur_{m-k})). \]
	By (\ref{eqn:eulerians}) and the inductive hypothesis, we have
	\[ \rk(\D_{\aaa}(\schur_m))=\sum_{b_{1}=1}^{a_{1}-1}\sum_{\bbb'\leq\aaa'}\E(\overleftarrow{\aaa'},\overleftarrow{\bbb'})\E(a_{1}-1,b_{1}-1)\sum_{k=1}^m \binom{k+b_{1}-1}{a_{1}-1}\binom{r+m-k+|\bbb'|-\nu_1-1}{r+|\aaa'|-1}, \]
	where we write $\bbb'=(b_2,\dots,b_h)$.
	As before, the last sum has non-trivial summands only for $a_{1}-b_{1}\leq k \leq m-|\aaa'|+|\bbb'|-\nu_1$ and we can apply Lemma~\ref{lemma:binomial for sym} with $A=a_{1}-1$, $R=r+|\aaa'|-1$ and $N=m-|\aaa|+|\bbb|-\nu_1$ (here $\bbb=(b_1,\dots,b_{h})$), resulting in 
	\[ \sum_{k=1}^m\binom{k+b_{1}-1}{a_{1}-1}\binom{r+m-k+|\bbb'|-\nu_1-1}{r+|\aaa'|-1}=\binom{r+m+|\bbb|-\nu_1-1}{r+|\aaa|-1}, \]
	which is again (\ref{eqn:in the proof of ch of sym}).
	
	Finally, if $a_1=\cdots=a_h=1$, then we can either proceed as above (using the formula for $\rk\D_1(\schur_m)$ above) or we can directly compute $\rk\D_1^{h}(\schur_m)$ thanks to Example~\ref{example:D 1 k}. In both cases we get
	\[ \rk\D_1^{h}(\schur_m)=\binom{r+m-1}{r+h-1}, \]
	which again agrees with (\ref{eqn:in the proof of ch of sym}), thus concluding the proof. \qed

	\subsection{Proof of Proposition~\ref{prop:ch of wedge}}
	
	We will need the following well-known binomial identity.
	
	\begin{lemma}\label{lemma:binomial for wedge}
		For any non-negative integers $A$, $R$ and $N$,
		\begin{equation*}
			\sum_{i=0}^N\binom{A+i}{A}\binom{R}{N-i}=\binom{R-A-1}{N}.
		\end{equation*}
	\end{lemma}
	
	The proof goes as in \S~\ref{section:proof of sym}, so we will only comment on the due changes. 
	First of all, as before it is enough to prove that
	\begin{equation}
		\rk(\D_{\aaa_\nu}(\schur_{1^m}))=\sum_{1\leq\bbb\leq\aaa_\nu}\E(\overleftarrow{\aaa_\nu},\overleftarrow{\bbb})\binom{r-|\aaa_\nu|}{m-|\aaa_\nu|+|\bbb|}.
	\end{equation}
	This already implies the formula for the Chern character, while the formula for the twisted Chern character is obtained by noticing that $k_1(E)=0$ as before.
	
	Let us start by settling the base of the induction and prove the case $\aaa_\nu=(a)$. As before, we have different behaviour if $a=1$ or $a>1$:
	\[ \rk\D_a\left(\extp^mE\right)=\sum_{k=1}^m (-1)^{k+1} k^{a-1}\binom{r}{m-k}=\left\{
	\begin{array}{ll}
		\binom{r-1}{m-1} & \mbox{if } a=1 \\
		\sum_{b=1}^{a-1}\E(a-1,b-1)\binom{r-a}{m-a-b} & \mbox{if } a>1
	\end{array}\right. \]
	where we used Lemma~\ref{lemma:binomial for wedge}.
	
	Put now $\aaa_\nu=(a_1,\dots,a_h)$ and assume that for at least one index $i$ we have $a_i>1$. Without loss of generality, it is $a_1>1$. If $\aaa'=(a_2,\dots,a_h)$, then
	\begin{align*} 
		\rk\D_{\aaa_\nu}\left(\extp^mE\right)= & \, \sum_{k=1}^m (-1)^{k+1} k^{a_1-1}\rk\D_{\aaa'}\left(\extp^{m-k}E\right) \\
		= & \, \sum_{b_1=1}^{a_1-1}\sum_{\bbb'\leq\aaa'}\E(\overleftarrow{\aaa'},\overleftarrow{\bbb'})\E(a_1-1,b_1-1)\sum_{k=1}^m (-1)^{k+1} \binom{k+b_1-1}{a_1-1}\binom{r-|\aaa'|}{m-|\aaa'|+|\bbb'|-\nu_1} \\
		= & \, \sum_{\bbb\leq\aaa_{\nu}}\E(\overleftarrow{\aaa_\nu},\overleftarrow{\bbb})\binom{r-|\aaa_\nu|}{m-|\aaa_\nu|+|\bbb|-\nu_1},
	\end{align*}
	where the last equality follows from Lemma~\ref{lemma:binomial for wedge}.
	
	Finally, the case when $a_1=\cdots=a_h=1$ is done similarly and it yields
	\[ \rk\D_1^h\left(\extp^mE\right)=\binom{r-1}{r-h}, \]
	thus concluding the proof. \qed
	

	\bibliographystyle{alpha}
	
\end{document}